\documentclass[12pt]{article}
\usepackage{amssymb,amsmath}
\begin{document}\thispagestyle{empty}

\begin{flushright}\noindent
Submitted October 1999; Revised December, 2001 \hfill\\ To appear
in {\it Compositio Mathematica}\hfill\\

\end{flushright}

\vspace*{2mm}

\begin{center}{\Large\bf
Resolution of Some Open Problems Concerning \\[3pt] Multiple Zeta
Evaluations of Arbitrary Depth }\\

\vglue 15mm {\bf Douglas Bowman}\footnote{Research partially
supported by NSF grant DMS-9705782.}
 \hglue 5mm{\tt
bowman}@{\tt math.uiuc.edu}\\\vglue 5mm University of Illinois at
Urbana-Champaign, Department of Mathematics,\\ 273 Altgeld Hall,
1409 W. Green St., Urbana, IL 61801 U.S.A.\\\vglue 3mm{\tt
http://www.math.uiuc.edu/$\!\sim$bowman}\\

\vglue 6mm

and

\vglue 6mm {\bf David M.~Bradley}\footnote{Research supported by
the University of Maine summer faculty research fund.} \hglue
5mm{\tt bradley}@{\tt math.umaine.edu }\\\vglue 5mm University of
Maine, Department of Mathematics and Statistics,\\ 5752 Neville
Hall, Orono, ME 04469-5752 U.S.A.\\\vglue 3mm{\tt
http://www.umemat.maine.edu/faculty/bradley/index.html}\\
\end{center}

\vglue10mm \noindent{\bf Abstract.} We prove some new evaluations
for multiple polylogarithms of arbitrary depth.   The simplest of
our results is a multiple zeta evaluation one order of complexity
beyond the well-known Broadhurst-Zagier formula.  Other results we
provide  settle three of the remaining outstanding conjectures of
Borwein, Bradley, and Broadhurst~\cite{BBB,BBBLa}. A complete
treatment of a certain arbitrary depth class of periodic
alternating unit Euler sums is also given.

\newcommand{\df}[2]{\mbox{$\frac{#1}{#2}$}}
\newcommand{\eu}{\stackrel{?}{=}}
\newcommand{\G}{\Gamma}
\newcommand{\g}{\gamma}
\newcommand{\Z}{{\mathbf Z}}
\newcommand{\R}{{\mathbf R}}
\newcommand{\Q}{{\mathbf Q}}
\newcommand{\C}{{\mathbf C}}
\newcommand{\w}{\omega}
\newcommand{\al}{\alpha}
\newcommand{\z}{\zeta}
\newcommand{\ou}{\overline{1}}
\newcommand{\tab}{\hskip.5cm}
\newcommand{\csch}{\mathrm{csch}}
\newcommand{\eop}{\hfill$\square$}

\newtheorem{Lem}{Lemma}
\newtheorem{Thm}{Theorem}
\newtheorem{Cor}{Corollary}
\newtheorem{Prop}{Proposition}

\section{Introduction}
\label{sect:motive}

The study of special values of multiple zeta functions concerns
itself with relations between values at integer vectors
$(s_1,\dots,s_k)$ of sums of the form
\begin{equation}
   \z(s_1,\dots,s_k):=\sum_{n_1>\cdots>n_k>0}
   \;\prod_{j=1}^k n_j^{-s_j},
\label{MZVdef}
\end{equation}
commonly referred to as {\it multiple zeta
values}~\cite{BBBLa,Hoff4,YOhno,Zag}.   We are primarily
interested in positive integer values of the arguments
$s_1,\dots,s_k$, in which case it is easily seen that $s_1>1$ is
both necessary and sufficient for the sum~(\ref{MZVdef}) to
converge.

A good deal of work on multiple zeta values has focused on the
problem of determining when ``complicated'' sums can be expressed
in terms of ``simpler'' sums. A crude but convenient measure of
the complexity of the sum~(\ref{MZVdef}) is the number $k$ of
nested summations.  This is also equal to the number of arguments
in the definition~(\ref{MZVdef}), and is called the {\it depth}.
Thus, researchers are interested in determining which sums can be
expressed in terms of other sums of lesser depth.

For each positive integer $k$, let $Z_k$ denote the set of all
multiple zeta values (with positive integer arguments) of depth
less than $k$. If we restrict our attention to rational polynomial
relationships, the aforementioned problem amounts to determining
for each $k$ the values of the arguments $s_1,\dots,s_k$ for which
the sum~(\ref{MZVdef}) lies in the polynomial ring $\Q[Z_k]$.
Settling this question in complete generality is currently well
beyond the reach of number theory.  For example, there is as yet
no proof that $\z(5)\notin\Q[Z_1]=\Q$, although it is strongly
suspected that $\z(5)$ is indeed irrational. Nevertheless,
considerable progress has been made with regard to proving
specific classes of reductions, even at arbitrary depth. A brief
historical overview will serve to put the problem in perspective.

Apart from Euler's celebrated depth-1 evaluation for the Riemann
zeta function
\[
   \z(2n) = \sum_{j=1}^\infty \frac1{j^{2n}}
   = -\frac12\cdot\frac{(2\pi i)^{2n}B_{2n}}{(2n)!},
   \qquad 0\le n\in \Z,
\]
in terms of the Bernoulli numbers $B_0=1,B_1=-1/2, B_2 = 1/6,
B_3=0, B_4=-1/30$, etc. defined by
\[
   \frac{z}{e^z-1}= \sum_{n=0}^\infty \frac{B_n}{n!} z^n,
   \qquad |z|<2\pi,
\]
the study of multiple zeta values can be fairly said to have begun
with Euler's depth-2 reduction~\cite{LE}
\begin{equation}
  2\z(n,1)=n\z(n+1)-\sum_{k=1}^{n-2}\z(n-k)\z(k+1),
   \qquad 2\le n\in\Z,
\label{EulerDepth2}
\end{equation}
expressing an infinite class of multiple zeta values of depth 2 in
terms of depth-1 (Riemann zeta) values. More generally, we refer
to a relation amongst multiple zeta values as a {\it depth-$k$
reduction} if the relation expresses a multiple zeta value of
depth $k$ in terms of multiple zeta values of depth less than $k$.

The first systematic study of reductions up to depth 3 was carried
out by Borwein, Bailey and Girgensohn, in a short series of
papers~\cite{Bail,BBG,BG} appearing in the early 1990s.  Research
undertaken by Hoffman~\cite{Hoff1,Hoff2,Hoff3,Hoff4},
Zagier~\cite{Zag}, and Borwein-Bradley-Broadhurst~\cite{BBB} on
multiple zeta values of arbitrary depth led to the discovery of
several relations satisfied by them. These relations can be
exploited by computer algebra systems to prove reductions of small
weight~\cite{MinPet1}.  (Here the {\it weight} of the multiple
zeta value~(\ref{MZVdef}) is simply the sum of the arguments:
$s_1+s_2+\cdots+s_k$.)

Additionally, high-precision evaluation of specific multiple zeta
values combined with human-directed computer searches using
lattice basis reduction algorithms led to several beautiful
conjectures concerning {\it arbitrary} depth
reductions~\cite{BBB,BBBLa,BBBLc}. An example of an arbitrary
depth reduction is
\begin{equation}
   \z(\underbrace{3,1,3,1,\dots,3,1}_{2n})
   =4^{-n}\z(\underbrace{4,4,\dots,4}_{n})
   = \frac{2^{1-4n}\z(4n)}{(2n+1)(4n+1)|B_{4n}|}
   = \frac{2\pi^{4n}}{(4n+2)!},
\label{Z31}
\end{equation}
in which the positive integers $2n$ and $n$ beneath the
underbraces in~(\ref{Z31}) denote the depth of the respective
multiple zeta values.  The formula~(\ref{Z31}), originally
conjectured by Zagier~\cite{Zag}, was first proved by Broadhurst.
A modification of Broadhurst's proof appears in~\cite{BBBLa}.
Subsequently, purely combinatorial proofs were
given~\cite{BBBLc,BowBradRyoo} based on the well-known shuffle
property of iterated integrals. To simplify the reading of such
formulas, when a string of arguments is repeated an exponent is
used.  In other words, we treat string multiplication as
concatenation. With this notation, the first two members of
formula~(\ref{Z31}) may be written
$\z(\{3,1\}^n)=4^{-n}\z(\{4\}^n).$

Although we focus here on analytic aspects of ultimately periodic
multiple polylogarithms, other aspects also play an important role
in related work.
See~\cite{Drin,Gonch,Gonch2,BK1,DJB1,BBBLc,BowBrad2,BowBradRyoo}
for connections with quantum groups, motivic Lie algebras,
algebraic geometry, knot theory, quantum field theory, and shuffle
combinatorics, respectively.  The survey~\cite{BowBrad3} provides
additional references and pointers to the literature.

\section{Statement of Results}
\label{sect:Main}

The Broadhurst-Zagier formula~(\ref{Z31}) is an example of an
arbitrary depth reduction in which the argument strings are
periodic. Whereas the case of periodic strings of period 1 is
quite well understood---there is a formula expressing
$\z(\{s\}^n)$ in terms of the depth-1 zeta values
$\z(s),\z(2s),\dots,\z(ns)$~\cite{BBB}---very little is known in
general about strings with longer periods apart from the fact that
their associated bivariate generating functions satisfy a
differential equation of order equal to the weight of the period.
In this paper, we offer some new multiple zeta evaluations for
ultimately periodic strings of period $2$ and period weight $4$.
We also give a complete treatment of a certain class of ultimately
periodic alternating unit Euler sums of period 2.  These results
are highlighted in subsection~\ref{sect:highlight} below. Three of
our results settle conjectures from~\cite{BBB} and~\cite{BBBLa}.

\subsection{Additional Notation}
\label{sect:notn}

Let $x$ be a real number satisfying $0\le x<1$.  The parametrized
multiple zeta function
\begin{equation}
   \z_x(s_1,\dots,s_k):=\sum_{n_1>\cdots>n_k>0}x^{n_1}
   \prod_{j=1}^k n_j^{-s_j}
\label{MLiDef}
\end{equation}
is defined for positive integers $s_j$, and is an instance of a
multiple polylogarithm~\cite{BBBLa,BowBrad3,Gonch}. Of course, if
$s_1>1$, then we can allow $x=1$ and in that case,~(\ref{MLiDef})
coincides with~(\ref{MZVdef}).  Euler sums have the form
\begin{equation}
   \z(s_1,\dots,s_k) := \sum_{n_1>\cdots>n_k>0}\; \prod_{j=1}^k
   n_j^{-|s_j|}\sigma_j^{-n_j},
\label{AltEuler}
\end{equation}
where $s_1,\dots,s_k$ are non-zero integers and
$\sigma_j:=\mathrm{signum}(s_j)$.  Thus, a multiple zeta value is
an Euler sum with no alternations, i.e.\ each $\sigma_j=1$.   If
each $s_j=\pm 1$, we also refer to~(\ref{AltEuler}) as a {\it
unit} Euler sum. To avoid confusion with the notion of analytic
continuation, we shall henceforth adopt the notation
of~\cite{BBB}, in which each $s_j$ in~(\ref{AltEuler}) is replaced
by $\overline{-s_j}$ when $s_j<0$. Thus, for example,
$\z(\ou)=-\log 2.$

Let $a,b,c$ and $x$ be complex numbers with $|x|<1$ and $c$ not
equal to zero or a negative integer.  We denote the Gaussian
hypergeometric function by
\begin{equation}
   F(a,b;c;x) = {}_2F_1\bigg(\begin{array}{cc}a,b\\c
   \end{array}\bigg|x\bigg) = \sum_{n=0}^\infty \frac{(a)_n
   (b)_n}{(c)_n}\cdot
     \frac{x^n}{n!},
\label{Fdef}
\end{equation}
where
\[
   (a)_n := \begin{cases} \displaystyle
   \prod_{j=0}^{n-1} (a+j), &\text{if\, $n$
   is a positive integer,}\\
   1, &\text{if\, $n=0$}\end{cases}
\]
is the Pochammer symbol (rising factorial).  If $\Re(c-b-a)>0$,
the series~(\ref{Fdef}) converges even when $|x|=1$, and thus for
complex $z$ we may abbreviate
\begin{align}
   Y_1(x,z) &:= F(z,-z;1;x),  &|x| \le 1,\label{Y1def}\\
   Y_2(x,z) &:= (1-x)F(1+z,1-z;2;1-x), &|1-x|<1.\label{Y2def}
\end{align}
Also, for complex $z$ not a positive integer multiple of $\pm 1$
or $\pm i$, we set
\begin{equation}
   G(z) :=
   \tfrac14\left\{\psi(1+iz)+\psi(1-iz)-\psi(1+z)-\psi(1-z)\right\},
   \label{Gdef}
\end{equation}
where $\psi=\Gamma'/\Gamma$ is the logarithmic derivative of the
Euler gamma function. Finally, for any non-negative integer $k$,
power series coefficient extraction will be performed throughout
by
\[
   [t^k]\sum_{n=0}^\infty a_n t^n = a_k.
\]
The notation will not be confused with the truncated square
brackets in $\lfloor x \rfloor$, which we use to denote the
greatest integer not exceeding the real number $x$.

\subsection{Main Results}
\label{sect:highlight} Our signal result is the following identity
for the generating function of the ultimately periodic sequence of
multiple polylogarithms $\{\z_x(3,\{1,3\}^n):0\le n\in\Z\}$.

\begin{Prop}\label{prop:zx313gf} For each real $x$ satisfying $0\le
x\le 1$, the formal power series
\begin{equation}
   S(x,z) := \sum_{n=0}^\infty (-1)^n z^{4n+2} 4^n
             \z_x(3,\{1,3\}^n)\label{z313GFdef}
\end{equation}
defines an entire function of the complex variable $z$.
Furthermore, if $Y_1$, $Y_2$ and $G$ are defined as
in~\textup{(\ref{Y1def}), (\ref{Y2def}) and (\ref{Gdef})}, then we
have the identity
\begin{equation}
   S(x,z)=G(z)Y_1(x,z)Y_1(x,iz)-\frac{Y_1(x,iz)Y_2(x,z)}{4Y_1(1,z)}
   +\frac{Y_1(x,z)Y_2(x,iz)}{4Y_1(1,iz)}
\label{z313GF}
\end{equation}
for all pairs $(x,z)$ for which the right hand side is defined.
\end{Prop}
Proposition~\ref{prop:zx313gf} has several interesting
consequences, two of which are given by Theorems~\ref{thm:z313}
and~\ref{thm:z213} below. Our first theorem gives a multiple zeta
value reduction which seems to have escaped the extensive
numerical and symbolic searches carried out by \cite{BBB} and
\cite{BBBLa}.

\begin{Thm}\label{thm:z313} For all non-negative integers $n$,
\begin{setlength}{\multlinegap}{60pt}
\begin{multline*}
   \z(3,\{1,3\}^n) = 4^{-n}\sum_{k=0}^n
   (-1)^k\z(4k+3)\z(\{4\}^{n-k})\\
   = \sum_{k=0}^n
   \frac{2\pi^{4k}}{(4k+2)!}\left(-\frac{1}{4}\right)^{n-k}\z(4n-4k+3).
\end{multline*}
\end{setlength}
\end{Thm}

\begin{Thm}[Conjectured in~\cite{BBB} and~\cite{BBBLa}]\label{thm:z213}
For all non-negative integers $n$,
\begin{eqnarray*}
   &&\z(2,\{1,3\}^n)\\
   &&\qquad = 4^{-n}\sum_{k=0}^n(-1)^k\z(\{4\}^{n-k})
   \bigg\{(4k+1)\z(4k+2)-4\sum_{j=1}^k\z(4j-1)\z(4k-4j+3)\bigg\}.
\end{eqnarray*}
\end{Thm}
We have also obtained cognate results for alternating unit Euler
sums~\cite{BBB,BBBLa} and multiple
polylogarithms~\cite{BBBLa,BowBrad3,Gonch} at one-half.   These
results are proved in section~\ref{sect:unit}, in which we settle
two conjectures from~\cite{BBB} and as a consequence obtain
reductions for the arbitrary depth alternating unit Euler sums
$\z(\{\ou,1\}^{n})$, $\z(\ou,\{1,\ou\}^{n})$,
$\z(\ou,\{\ou,1\}^{n})$, and $\z(\ou,\ou,\{1,\ou\}^{n})$, where
$n$ is any non-negative integer. The main result from which these
reductions are derived is Proposition~\ref{prop:Mgf} below.

Recall the generating function
\begin{equation}
   A(z):=\sum_{n=0}^{\infty}z^n\z(\{\ou\}^n) =
   \prod_{j=1}^\infty \bigg(1+\frac{(-1)^jz}{j}\bigg)
   =\frac{\G(1/2)}{\G(1+z/2)\G(1/2-z/2)}
\label{Adef}
\end{equation}
from~\cite{BBB}.
\begin{Prop}\label{prop:Mgf}  For each real $x$ satisfying $0\le
x\le 1$, the formal power series
\begin{equation}
   M(x,t) :=  \sum_{n=0}^\infty
    \bigg[t^{2n}\z_x(\{\ou,1\}^n)+t^{2n+1}\z_x(\ou,\{1,\ou\}^n)\bigg]
\label{Mdef}
\end{equation}
defines an entire function of the complex variable $t$.
Furthermore, if we put $z=(1+i)t/2$, $s=(1+x)/2$, and let $U(s,z)
= Y_1(s,z)-zY_2(s,z)$, where $Y_1$ and $Y_2$ are given
by~\textup{(\ref{Y1def})} and~\textup{(\ref{Y2def})} respectively,
then
\begin{equation}
   M(x,t) = \frac{U(s,-z)U(s,iz)}{A(-z)A(iz)}
\label{Mrep}
\end{equation}
for all pairs $(x,t)$ for which the right hand side is defined.
\end{Prop}
From Proposition~\ref{prop:Mgf} we obtain the following generating
function identities for alternating unit Euler sums.

\begin{Cor}[Conjectured in~\cite{BBB}]\label{BBB14}
Let $A(z)$ be as in~\textup{(\ref{Adef})}.  Then for all complex
numbers $t$, we have
\[
   \sum_{n=0}^\infty
   \bigg[t^{2n}\z(\{\ou,1\}^n)+t^{2n+1}\z(\ou,\{1,\ou\}^n)\bigg]
   = A\left(\frac{t}{1-i}\right)A\left(\frac{t}{1+i}\right).
\]
\end{Cor}

\begin{Cor}[Conjectured in~\cite{BBB}]\label{Tmilk}
Let the functions $G$ and $A$ be given by~\textup{(\ref{Gdef})}
and~\textup{(\ref{Adef})}, respectively.  Then
\begin{equation}
\begin{split}
&1+\sum_{n=0}^\infty\bigg[
   t^{2n+1}\z(\ou,\{\ou,1\}^n)+t^{2n+2}\z(\ou,\ou,\{1,\ou\}^n)\bigg]\\
&\qquad\qquad=\tfrac12(1+i)zA(z)A(-iz)\big\{\pi\csc(\pi z)-i\pi\,
   \csch(\pi z)+4G(z)\big\},
\end{split}
\label{bbb16mod}
\end{equation}
holds for all complex numbers $t$ and $z$ satisfying $z=(1+i)t/2$
and such that the right hand side is defined.
\end{Cor}
\noindent{\bf Remark.} It is routine to check that the right hand
side of~(\ref{bbb16mod}) and the corresponding expression in
equation~(16) of~\cite{BBB} are equal. Our expression has the
advantage of being easier to express as a Maclaurin series.

From Corollaries~\ref{BBB14} and~\ref{Tmilk} we obtain the
following reductions for alternating unit Euler sums.

\begin{Thm}
\label{thm:cream} Define a sequence of numbers
$c_0,c_1,\ldots\in\Q[\log 2,\z(2),\z(3),\dots]$ by
\[
   \sum_{m=0}^\infty c_m\, x^m = \exp\bigg(\sum_{k=1}^\infty b_k
   \frac{x^k}{k}\bigg),
\]
where $b_1=-\log 2$ and
\[
    b_k=
\begin{cases}
    (-1)^{\lfloor(k+1)/4\rfloor}\, 2^{(1-k)/2}(2^{1-k}-1)\z(k),
     &\text{if\, $1<k$ is odd,}\\
    (-1)^{1+k/4}\, 2^{1-k/2}\z(k), &\text{if\, $k\equiv 0\bmod 4$,}\\
    0, &\text{if\, $k\equiv 2\bmod 4$.}
\end{cases}
\]
Then for all non-negative integers $n$ we have
\[
   \z(\{\ou,1\}^n)=c_{2n}, \quad \text{and} \quad
   \z(\ou,\{1,\ou\}^n)=c_{2n+1}.
\]
\end{Thm}

\begin{Thm}
\label{thm:cheese} Let $b_k$ be as in Theorem~\ref{thm:cream}.
Define numbers $c_m'\in\Q[\log 2,\z(2),\z(3),\dots]$ by
\[
   \sum_{m=0}^\infty c_m'\, x^m = \bigg(\sum_{k=0}^\infty d_k\,
   x^k\bigg)\exp\bigg(\sum_{k=1}^\infty b_k\frac{x^k}{k}\bigg),
\]
where $d_0=1$ and
\[
    d_k=
\begin{cases} (-1)^{\lfloor k+1/2\rfloor}\, 2^{2-3k/2}(2^{k-1}-1)\z(k),
   &\text{if\, $0<k$ is even,}\\
    0, &\text{if\, $k\equiv 1\bmod 4$,}\\
    (-1)^{(k+1)/4}\,2^{(3-k)/2}\z(k), &\text{if\, $k\equiv 3\bmod 4$.}
\end{cases}
\]
Then for all non-negative integers $n$ we have
\[
   \z(\ou,\{\ou,1\}^n)=c_{2n+1}',
   \quad \text{and}\quad
   \z(\ou,\ou,\{1,\ou\}^n)=c_{2n+2}'.
\]
\end{Thm}

\noindent{\bf Remark.}  Of course, one can express the numbers
$c_m$ and $c_m'$ of Theorems~\ref{thm:cream} and~\ref{thm:cheese}
explicitly. Thus,
\[
  c_m' = \sum_{k=0}^m c_k d_{m-k} \quad \text{and}\quad
  c_m
      = \sum\prod_{k\ge
      1}\frac{1}{j_k!}\bigg(\frac{b_k}{k}\bigg)^{j_k},
\]
where the sum is over all non-negative integers $j_1,j_2,\dots$
satisfying $\sum_{k\ge 1}kj_k=m$.

\section{Multiple Zeta Values of Period 2}
\label{sect:z213}

This section contains the proofs of our results pertaining to
ultimately periodic multiple zeta values of period 2 and period
weight 4, namely Proposition~\ref{prop:zx313gf} and
Theorems~\ref{thm:z313} and~\ref{thm:z213}.  It is interesting to
note that the proof of Propositions~\ref{prop:zx313gf}
and~\ref{prop:Mgf} relies on a general result which also
eliminates the need for computer algebra in the orginal proof of
the Broadhurst-Zagier formula in~\cite{BBBLa}.  See
Lemma~\ref{lem:zfac} below.

\subsection{Proof of
Proposition~\ref{prop:zx313gf}}\label{sect:zfac} Let $R(x,z)$
denote the right hand side of~(\ref{z313GF}).  For each fixed $x$
satisfying $0<x\le 1$, $R(x,z)$ is evidently an analytic function
of $z$, apart from isolated singularities (possible poles) at the
positive integer multiples of $\pm 1$, $\pm i$.    Let $n$ be a
positive integer and let $p\in\{1,-1,i,-i\}$.  A straightforward
calculation shows that $\lim_{z\to pn}(z-pn)R(x,z) = 0$ follows
from the identity
\begin{equation}
   Y_1(x,n)+n(-1)^n Y_2(x,n)=0,
   \qquad 1\le n\in\Z,
\label{JacobiY}
\end{equation}
which in turn is a consequence of the identity~\cite[p.\
254]{Rain}
\begin{setlength}{\multlinegap}{30pt}
\begin{multline*}
   \frac{(1+\alpha)_n}{n!}\;
   {}_2F_1\bigg(\begin{array}{cc}-n,1+\alpha+\beta+n\\1+\alpha
   \end{array}\bigg|\frac{1-y}{2}\bigg)\\
   = \frac{(-1)^n(1+\beta)_n}{n!}\;
   {}_2F_1\bigg(\begin{array}{cc}-n,1+\alpha+\beta+n\\1+\beta\end{array}
   \bigg|\frac{1+y}{2}\bigg)
\end{multline*}
\end{setlength}for the Jacobi polynomials.  Thus, for $0<x\le 1$, the
singularities of $R(x,z)$ are all removable.

It now suffices to show that $S(x,z)$ and $R(x,z)$ both have
Maclaurin series in $x$ which begin
\[
   z^2x+\tfrac18z^2x^2+\tfrac{1}{27}z^2\big(1-2z^4\big)x^3
   +\tfrac{1}{64}z^2\big(1-\tfrac72z^4\big)x^4+O(x^5),
   \qquad x\to0,
\]
and are both annihilated by the  differential operator
\begin{equation}
   D_1^2D_0^2 +4z^4,\qquad
   D_0:=x\frac{d }{d x},\quad D_1:=(1-x)\frac{d
}{d x}. \label{dop}
\end{equation}
Checking these facts for $S(x,z)$ is a trivial exercise.  Although
in general $R(x,z)$ is undefined when $x=0$ because
$Y_2(0,z)=F(1+z,1-z;2;1)$ diverges unless it terminates, a short
calculation employing the special case
\begin{multline*}
   F(1+z,1-z;2;x) \\= \frac{\Gamma(2)}{\Gamma(1+z)\Gamma(1-z)}
     \sum_{n=0}^\infty\frac{(1+z)_n(1-z)_n}{(n!)^2}
   \bigg\{2\psi(n+1)-\psi(n+1-z)\\
   -\psi(n+1+z)-\log(1-x)\bigg\}(1-x)^n,
   \qquad 0\le x<1
\end{multline*}
of the formula 15.3.10 of~\cite[p.\ 559]{AS} verifies that
$R(x,z)$ possesses a Maclaurin series at $x=0$ that begins as
stated.

The fact that $(D_1^2D_0^2 +4z^4)R(x,z)=0$ is most easily seen by
setting $f(x)=1-x$, $g(x)=x$, and $t=z^2$ in Lemma~\ref{lem:zfac}
below. Two solutions to the
 differential equation $(D_1D_0+z^2)y=0$ are given by
$y=Y_1(x,z)$ and $y=Y_2(x,z)$, and thus changing $z$ to $iz$ we
see that two solutions of the  differential equation
$(D_1D_0-z^2)y=0$ are $y=Y_1(x,iz)$ and $y=Y_2(x,iz)$. The lemma
then shows that each of the functions $Y_1(x,z)Y_1(x,iz)$,
$Y_1(x,iz)Y_2(x,z)$, and $Y_1(x,z)Y_2(x,iz)$ are annihilated by
the operator~(\ref{dop}). \eop

\begin{Lem}\label{lem:zfac} Let $K$ be a field of
characteristic not equal to $2$ and let $D$ be a derivation on
$K$. For each $k\in K,$ define a derivation $D_k:=kD$.  Let $t$ be
a constant, and suppose that for some $f,g,u,v\in K$ the
differential equations $(D_fD_g+t)u=0$ and $(D_fD_g-t)v=0$ hold.
Then $uv$ is annhilated by the differential operator
$(D_f^2D_g^2+4t^2).$
\end{Lem}
\noindent{\bf Proof of Lemma~\ref{lem:zfac}.} First note that
$uD_g^2v+vD_g^2u=0$, for
\[
   \left(uD_g^2v+vD_g^2u\right)f
  =\left(uD_fD_gv+vD_fD_gu\right)g
  =(utv-vtu)g=0.
\]
We now calculate $D_f^2D_g^2(uv)$. By the Leibniz rule and our
note above,
\[
   D_f^2D_g^2(uv)
  =D_f^2\left(uD_g^2v+2(D_gu)(D_gv)+vD_g^2u\right)
  =2D_f^2(D_gu)(D_gv).
\]
But
\begin{align*}
   D_f^2(D_gu)(D_gv)
   &= D_f(D_gu)(D_fD_gv) + D_f(D_fD_gu)(D_gv)\\
   &= D_f(D_gu)(tv) + D_f(-tu)(D_gv)\\
   &= t\left[vD_fD_gu+(D_fv)(D_gu) -(D_fu)(D_gv)-
uD_fD_gv\right].
\end{align*}
In the previous expression, the middle two terms cancel since
\[
   (D_fv)(D_gu)=(fg)(Dv)(Du)=(D_fu)(D_gv).
\]
Hence we have
\[
  D_f^2D_g^2(uv)
  = 2t\left[vD_fD_gu-uD_fD_gv\right]
  = 2t\left[v(-tu)-u(tv)\right]=-4t^2uv.
\]
\eop

\subsection{Proof of Theorem~\ref{thm:z313}}
Let $x=1$ in the identity~(\ref{z313GF}) of
Proposition~\ref{prop:zx313gf} and then extract the coefficient of
the appropriate power of $z$. More explicitly, we note that
\begin{equation}
   \z(3,\{1,3\}^n)=(-1)^n 4^{-n} [z^{4n+2}]S(1,z).
\label{z313extract}
\end{equation}
But
\[
   S(1,z) = G(z)Y_1(1,z)Y_1(1,iz) = G(z)Q(z),
\]
where by Gauss's ${}_2F_1$ summation theorem and the infinite
product formula for sine,
\begin{equation}
   Q(z) := Y_1(1,z)Y_1(1,iz)
   =\frac{\sin(\pi z)}{\pi z}\cdot\frac{\sin(\pi i z)}{\pi i z}
   = \sum_{n=0}^\infty (-1)^n z^{4n}\z(\{4\}^n).
\label{Qdef}
\end{equation}
Since
\begin{equation}
   G(z) = \sum_{n=0}^\infty z^{4n+2}\z(4n+3),\qquad |z|<1,
\label{Gexpand}
\end{equation}
the desired formula follows from~(\ref{z313extract}),~(\ref{Z31}),
and the formula for the coefficient in the Cauchy product of power
series. \eop

\begin{Cor}\label{cor:z2131} For all non-negative integers $n$,
\begin{setlength}{\multlinegap}{60pt}
\begin{multline*}
  \z(2,1,\{3,1\}^n) =  4^{-n}\sum_{k=0}^n
  (-1)^k\z(4k+3)\z(\{4\}^{n-k})\\
   = \sum_{k=0}^n
   \frac{2\pi^{4k}}{(4k+2)!}\left(-\frac{1}{4}\right)^{n-k}\z(4n-4k+3).
\end{multline*}
\end{setlength}
\end{Cor}

\noindent{\bf Proof.} Apply duality~\cite{BBBLa,Kass,Zag} to
Theorem~\ref{thm:z313}. \eop

\subsection{Proof of Theorem~\ref{thm:z213}}
Differentiate both sides of~(\ref{z313GF}) in
Proposition~\ref{prop:zx313gf} with respect to $x$, let $x\to1-$
(requires asymptotic formulas for the relevant hypergeometrics)
and then extract the coefficient of the appropriate power of $z$.
More explicitly, we note that
\[
   x\frac{d}{d x}S(x,z)
   = \sum_{n=0}^\infty (-1)^n z^{4n+2} 4^n\z_x(2,\{1,3\}^n),
\]
and hence (letting prime denote differentiation with respect to
the first argument)
\begin{equation}
   \z(2,\{1,3\}^n) = (-1)^n 4^{-n}[z^{4n+2}]S'(1,z).
\label{z213extract}
\end{equation}
Differentiating~(\ref{z313GF}) we get that
\begin{setlength}{\multlinegap}{30pt}
\begin{multline}
   S'(x,z) = G(z)\{H(x,z)+H(x,iz)\}
   -\frac{Y_1'(x,iz)Y_2(x,z)}{4Y_1(1,z)}\\
   -\frac{Y_1(x,iz)Y_2'(x,z)}{4Y_1(1,z)}
   +\frac{Y_1'(x,z)Y_2(x,iz)}{4Y_1(1,iz)}
   +\frac{Y_1(x,z)Y_2'(x,iz)}{4Y_1(1,iz)},
\label{Sprime}
\end{multline}
\end{setlength}where
\begin{equation}
   H(x,z) := -z^2F(1+z,1-z;2;x)F(iz,-iz;1;x).
\label{Hdef}
\end{equation}
Entries [15.3.10] and [15.3.11] of~\cite{AS} provide the
asymptotic formulas
\begin{align}
   F(1+z,1-z;2;x)&=\frac{2\psi(1)-\psi(1+z)-\psi(1-z)
   -\log(1-x)} {\G(1+z)\G(1-z)}\nonumber\\
   &+ O((1-x)\log(1-x)), \qquad x\to1-,
\label{Hfac1}
\end{align}
and
\begin{equation}
   F(iz,-iz;1;x)=\frac{1}{\G(1+iz)\G(1-iz)}+O((1-x)\log(1-x)),
   \qquad x\to1-,
\label{Hfac2}
\end{equation}
respectively.   If we now substitute the asymptotic
formulas~(\ref{Hfac1}) and~(\ref{Hfac2}) into~(\ref{Hdef}), apply
the reflection formula for the gamma function and the
definition~(\ref{Qdef}), there comes
\begin{align*}
   H(x,z) &= -z^2\frac{\sin(\pi z)}{\pi z}\cdot
            \frac{\sin(\pi i z)}{\pi i z}
   \bigg\{2\psi(1)-\psi(1+z)-\psi(1-z)-\log(1-x)+o(1)\bigg\}\\
   &= -z^2
   Q(z)\bigg\{2\psi(1)-\psi(1+z)-\psi(1-z)-\log(1-x)+o(1)\bigg\},
\end{align*}
and hence as $x\to1-$,
\begin{align}
   H(x,z)+H(x,iz) &=
   Q(z)\left\{z^2(\log(1-x)-2\psi(1)+\psi(1+z)+\psi(1-
z))\right.\nonumber\\
   &\left.\quad-z^2(\log(1-x)-2\psi(1)-\psi(1+iz)-\psi(1-
iz))\right\}\nonumber\\
   &\quad+o(1)\nonumber\\
   &= -4z^2Q(z)G(z)+o(1).
\label{Hsum}
\end{align}
We now substitute~(\ref{Hsum}) into~(\ref{Sprime}). Since
$Y_2(x,z)=O(1-x)$ and $Y_1'(x,z)=O(\log(1-x))$ as $x\to1-$, it
follows that
\begin{align}
   S'(x,z) &= -4z^2G(z)Q(z)G(z)\nonumber\\
    &-\frac{Y_1(x,iz)}{4Y_1(1,z)}
   \bigg\{-F(1+z,1-z;2;1-x)+O(1-x)\bigg\}\nonumber\\
   &+\frac{Y_1(x,z)}{4Y_1(1,iz)}\bigg\{-F(1+iz,1-iz;
   2;1-x)+O(1-x)\bigg\}+o(1).
\label{Slimit}
\end{align}
We now let $x\to1-$ in~(\ref{Slimit}), obtaining
\begin{equation}
   S'(1,z) = -4z^2Q(z)G^2(z)+\tfrac14\pi^2z^2Q(z)\left\{\csc^2(\pi
   z)-{\csch}^2(\pi z)\right\}.
\label{csch}
\end{equation}
In view of~(\ref{z213extract}), the proof of
Theorem~\ref{thm:z213} now follows on extracting the coefficient
of the appropriate power of $z$ from both sides of~(\ref{csch}).
\eop

\medskip\noindent{\bf Remark.} It is possible to continue
differentiating~(\ref{Sprime}) and obtain generating functions for
the multiple polylogarithms $\z_x(1,\{1,3\}^{n})$ and
$\z_x(\{1,3\}^{n})$. One can similarly differentiate
$Y_1(x,z)Y_1(x,iz)$ and obtain generating functions for
$\z_x(2,1,\{3,1\}^{n})$, $\z_x(1,1,\{3,1\}^{n})$, and
$\z_x(1,\{3,1\}^{n})$. Setting $x=1/2$ in all these generating
functions and dualizing then allows one to obtain reductions for
the alternating unit Euler sums $\z(\{\ou,1\}^{n})$,
$\z(\ou,\{1,\ou\}^{n})$, $\z(\ou,\{\ou,1\}^{n})$, and
$\z(\ou,\ou,\{1,\ou\}^{n})$ for all non-negative integers $n$. The
key step is to observe that the  derivatives
\[
   \left(\frac{d}{d x}\right)^k
   F(z,-z;1;x)\bigg|_{x=1/2}\qquad (0\le k\in \Z)
\]
can be expressed, via entry [15.1.25] of~\cite{AS}, in terms of
gamma functions, and thereby in terms of the generating
function~(\ref{Adef}) for the sequence $\{\z(\{\ou\}^n): 0\le n\in
\Z\}$, each term of which is reducible~\cite{BBB}. This procedure
is cumbersome, however, and in the next section we obtain the same
results by a much more elegant method. In the process we settle
two additional conjectures from~\cite{BBB}.

\section{Alternating Unit Euler Sums of Period 2}
\label{sect:unit} This section contains the proofs of our results
pertaining to ultimately periodic alternating unit Euler sums of
period 2, namely Proposition~\ref{prop:Mgf} and
Theorems~\ref{thm:cream} and~\ref{thm:cheese}.

\subsection{Proof of Proposition~\ref{prop:Mgf}}
One first checks that $M(x,t)$ satisfies the  differential
equation
\begin{equation}
   \left[(1-x)\frac{d}{d x}\right]^2
   \left[-(1+x)\frac{d}{d x}\right]^2
   M(x,t)=t^4M(x,t).
\label{Mde}
\end{equation}
To solve the  differential equation~(\ref{Mde}), one could apply
Lemma~\ref{lem:zfac} directly with $f(x)=1+x$ and $g(x)=1-x$.
However, it is more convenient to make a change of variable.  With
$s=(1+x)/2$, $z=(1+i)t/2$ and $L(s,z):=M(x,t)$,~(\ref{Mde}) goes
over into
\begin{equation}
   \left\{\left[(1-s)\frac{d}{d s}\right]^2
   \left[s\frac{d}{d s}\right]^2
   +4z^4\right\}L(s,z)=0,
\label{Lde}
\end{equation}
which we've already encountered~(\ref{dop}).  A routine
computation using entry [15.1.25] of~\cite{AS}
shows that
\begin{eqnarray*}
   U(1/2,z) &=& A(z),\\
   \frac{d}{d s}U(s,z)\bigg|_{s=1/2}&=&2zA(z),\\
   \frac{d^2}{d s^2}U(s,z)\bigg|_{s=1/2}&=&-4z(1+z)A(z),\\
   \frac{d^3}{d s^3}U(s,z)\bigg|_{s=1/2}&=&8z(1+z)(2-z)A(z).
\end{eqnarray*}
Using these relations, it is can be easily shown that the
Wronskian determinant of the four functions
\[
   U(s,z)U(s,iz),\quad U(s,-z)U(s,iz),\quad
   U(s,z)U(s,-iz),\quad U(s,-z)U(s,-iz)
\]
at $s=1/2$ is equal to
\[
   -2^{13}z^6\left(\frac{\sin(\pi z)}{\pi
   z}\right)^2\left(\frac{\sinh(\pi z)}{\pi z}\right)^2.
\]
Since the Wronskian is not identically zero, the four functions
are linearly independent. From the  differential
equation~(\ref{Lde}) and Lemma~\ref{lem:zfac} with $f(s)=1-s$ and
$g(s)=s$, it follows that there exist functions $\alpha(z)$,
$\beta(z)$, $\gamma(z)$, $\delta(z)$ such that
\begin{align}
   M(x,t)=L(s,z)&=\alpha(z)U(s,z)U(s,iz)+\beta(z)U(s,-z)U(s,iz)
   \nonumber\\
   &+\gamma(z)U(s,z)U(s,-iz)+\delta(z)U(s,-z)U(s,-iz).
\label{linearcombo}
\end{align}
Now setting $x=0$ (which is the same as $s=1/2$)
in~(\ref{linearcombo}), and performing the operations
\begin{eqnarray*}
   -(1+x)\frac{d }{d x}M(x,t)\bigg|_{x=0} &=&
   -s\frac{d }{d s}L(s,t)\bigg|_{s=1/2},\\
   \left[-(1+x)\frac{d }{d
   x}\right]^2M(x,t)\bigg|_{x=0}
   &=&\left[-s\frac{d }{d
   s}\right]^2L(s,z)\bigg|_{s=1/2},\\
   (1-x)\frac{d }{d x}\left[-(1+x)\frac{d }{d
x}\right]^2M(x,t)\bigg|_{x=0}
   &=&(1-s)\frac{d }{d s}\left[-s\frac{d }{d
s}\right]^2L(s,z)\bigg|_{s=1/2}
\end{eqnarray*}
yields the following system of equations:
\[
   \left[\begin{array}{rrrr}1& \qquad 1&\qquad 1&\qquad 1 \\
    -i&\qquad 1&\qquad -1&\qquad i\\-1&\qquad 1&\qquad 1&\qquad
-1 \\
    i&\qquad 1&\qquad -1&\qquad i\end{array}\right]
   \left[\begin{array}{l}A(z)A(iz)\alpha(z)\\A(-z)A(iz)\beta(z)\\
      A(z)A(-iz)\gamma(z)\\A(-z)A(-iz)\delta(z)\end{array}\right]
   =\left[\begin{array}{c}1\\1\\1\\1\end{array}\right].
\]
Inspection via Cramer's rule gives
$\alpha(z)=\gamma(z)=\delta(z)=0$ and
\[
   \beta(z)=\frac{1}{A(-z)A(iz)}.
\]
It follows that
\[
   M(x,t)=L(s,z)=\frac{U(s,-z)U(s,iz)}{A(-z)A(iz)},
\]
wherever the right hand side is defined.  As in the proof of
Proposition~\ref{prop:zx313gf}, we find that as a consequence of
the Jacobi polynomial identity~(\ref{JacobiY}), the singularities
(which in this case occur when $z$ is an even positive integer
multiple of 1 or $i$ or an odd positive integer multiple of $-1$
or $-i$) are all removable. Thus, $M(x,t)$ is an entire function
of $t$ for each $x$ satisfying $0\le x\le 1$. \eop

\subsubsection{Proof of Corollary~\ref{BBB14}} Set $x=1$ in
equation~(\ref{Mrep}) of Proposition~\ref{prop:Mgf}. In view of
the fact that $Y_2(1,z)=0$, we have $U(1,z)=Y_1(1,z)=\sin(\pi
z)/(\pi z)$. Thus,
\begin{align*}
   M(1,t) &= \frac{U(1,-z)U(1,iz)}{A(-z)A(iz)}\\
          &= \frac{Y_1(1,z)Y_1(1,iz)}{A(-z)A(iz)}\\
          &=
          A(z)A(-iz)\cdot\frac{Y_1(1,z)Y_1(1,iz)}{A(z)A(-z)A(iz)A(-iz)}\\
          &= A(z)A(-iz).
\end{align*}
\eop

\medskip\noindent{\bf Remark.} Broadhurst~\cite{Bprivate}
has outlined a different proof of Corollary~\ref{BBB14} using
iterated integrals.

\subsection{\bf Proof of Theorem~\ref{thm:cream}} By
Corollary~\ref{BBB14} we need to compute the Maclaurin series for
$$ A\left(\frac{t}{1-i}\right)A\left(\frac{t}{1+i}\right). $$ From
equation~(12) of~\cite{BBB},
\[
   A(z) = \exp\bigg(\sum_{k=1}^\infty\frac{(-1)^{k+1} a_k
z^k}{k}\bigg),
\]
where
\[
  a_k = {\mathrm{Li}_k}((-1)^k)
      = \sum_{n=1}^\infty\frac{(-1)^{nk}}{n^k}
      = \begin{cases}-\log 2, &\text{if $k=1$;}\\
                \z(k), &\text{if $k$ is even;}\\
          \left(2^{1-k}-1\right)\z(k), &\text{if $k>1$ is odd.}
          \end{cases}
\]
Hence
\[
   A\left(\frac{t}{1-i}\right)A\left(\frac{t}{1+i}\right)
   =\exp\bigg(\sum_{k=1}^\infty\frac{[(-1)^{k+1}-i^k]a_k t^k}
   {k(1-i)^k}\bigg).
\]
Simplifying the complex values in the sum we find that
\[
    A\left(\frac{t}{1-i}\right)A\left(\frac{t}{1+i}\right)
=\exp\bigg(\sum_{k=1}^\infty {b_k \frac{t^k}{k}}\bigg),
\]
where the sequence $b_k$ is as defined in the theorem. Expanding
the exponential completes the proof. \eop

Applying duality~\cite{BBBLa,Kass,Zag} to Theorem~\ref{thm:cream}
shows that we have also obtained reductions for certain multiple
polylogarithmic values at $1/2$. Specifically, for all
non-negative integers $n$,
\begin{eqnarray*}
   \z_{1/2}(\{3,1\}^n) &=& \z(\{\ou,1\}^{2n}),\\
   \z_{1/2}(2,1,\{3,1\}^n) &=& \z(\ou,\{1,\ou\}^{2n+1}),\\
   \z_{1/2}(1,1,\{3,1\}^n) &=& \z(\{\ou,1\}^{2n+1}),\\
   \z_{1/2}(1,\{3,1\}^n) &=& -\z(\ou,\{1,\ou\}^{2n}).
\end{eqnarray*}

\subsection{Proof of Corollary~\ref{Tmilk}}

For each real $x$ satisfying $0\le x\le 1$, form the formal power
series
\begin{equation}
   T(x,t) := 1+\sum_{n=0}^\infty
    \bigg[ t^{2n+1}\z_x(\ou,\{\ou,1\}^n)
          +t^{2n+2}\z_x(\ou,\ou,\{1,\ou\}^n)\bigg].
   \label{Tdef}
\end{equation}
Then it is a routine computation to show that
\[
   T(x,t)=-\bigg(\frac{1+x}{t}\bigg)\frac{d}{d x} M(x,t),
\]
where $M(x,t)$ is as in equation~(\ref{Mdef}).  Put $z=(1+i)t/2$,
$s=(1+x)/2$, and $U(s,z) = Y_1(s,z)-zY_2(s,z)$, where $Y_1$ and
$Y_2$ are given by~(\ref{Y1def}) and~(\ref{Y2def}), respectively.
From Proposition~\ref{prop:Mgf}, it follows that
\begin{equation}
\label{Tgf}
   T(x,t)=-\frac{s}{z}\bigg(\frac{1+i}{2}\bigg)
   \frac{U^{'}(s,iz)U(s,-z)+U(s,iz)U^{'}(s,-z)}{A(-z)A(iz)},
\end{equation}
wherever the right hand side is defined, and where the prime
denotes differentiation with respect to $s$.  As in the proof of
Proposition~\ref{prop:Mgf}, we find that for each $x$ satisfying
$0\le x\le 1$, $T(x,t)$ defines an entire function of $t$.

We now compute the right hand side of~(\ref{Tgf}) when $x=s=1$.
The obstacle to overcome is the singularity of $Y_1^{'}(s,az)$ at
$s=1$. (We will be taking $a$ to be $-1$ or $i$ as needed.) Of
course this function occurs twice in the generating function so
that the singularities cancel. In particular, from the definitions
of the functions it is immediate that
\[
  U^{'}(s,az)=az-a^2z^2F(1+az,1-az;2;s)+O(1-s),\qquad s\to 1-.
\]
By entry [15.3.10] of~\cite{AS} it follows that
\begin{eqnarray*}
    U^{'}(s,az)&=&az-\frac{(a^2z^2)(2\psi(1)-\psi(1+az)-\psi(1-az)
   -\log(1-s))}{\G(1+az)\G(1-az)}\\ &+& O((1-s)\log(1-s)),
   \qquad s\to1-.
\end{eqnarray*}
Similarly, from [15.3.11] of~\cite{AS},
\[
   U(s,az)=\frac{1}{\G(1+az)\G(1-az)}+ O((1-s)\log(1-s)),
   \qquad s\to 1-.
\]
Substituting these expressions (with $a$ taken to be $-1$ or $i$
as appropriate) into~(\ref{Tgf}), applying the reflection formula
for the gamma function, and simplifying yields
\begin{setlength}{\multlinegap}{40pt}
\begin{multline*}
   T(x,t) =\tfrac12(1+i)szA(z)A(-iz)\big\{\pi\csc(\pi z)-i\pi\,
   \csch(\pi z)+4G(z)\big\}\\+O((1-s)\log(1-s)),
   \qquad s\to1-.
\end{multline*}
\end{setlength}
Letting $s\to 1-$ completes the proof of the corollary.\eop

\subsection{Proof of Theorem~\ref{thm:cheese}} We now conclude
the proof of Theorem~\ref{thm:cheese}. Changing variables from $z$
to $t$ in Corollary~\ref{Tmilk}, using~(\ref{Gexpand}) and the
well-known Maclaurin series for cosecant gives
\begin{multline*}
   1+\sum_{n=0}^\infty\bigg[
   t^{2n+1}\z(\ou,\{\ou,1\}^n)+t^{2n+2}\z(\ou,\ou,\{1,\ou\}^n)\bigg]
   =A\left(\frac{t}{1-i}\right)A\left(\frac{t}{1+i}\right)\\
      \times\sum_{n=0}^\infty \bigg[(-1)^{\lfloor (n+1)/2\rfloor}
      2^{2-3n}(2^{2n-1}-1)\z(2n)t^{2n}
   -(-4)^{-n}\z(4n+3)t^{4n+3}\bigg].
\end{multline*}
But this last expression is exactly
$
   \sum_{k=0}^\infty c_kt^k\sum_{j=0}^\infty d_j t^j.
$
The theorem now follows by taking the Cauchy product and equating
coefficients.
\eop

Duality~\cite{BBBLa,Kass,Zag} applied to Theorem~\ref{thm:cheese}
provides evaluations for multiple polylogarithms at $1/2$.    For
all non-negative integers $n$,
\begin{eqnarray*}
   \z_{1/2}(3,\{1,3\}^n) &=& -\z(\ou,\{\ou,1\}^{2n+1}),\\
   \z_{1/2}(2,\{1,3\}^n) &=& -\z(\ou,\ou,\{1,\ou\}^{2n}),\\
   \z_{1/2}(1,\{1,3\}^n) &=& -\z(\ou,\{\ou,1\}^{2n}),\\
   \z_{1/2}(\{1,3\}^{n+1}) &=& \z(\ou,\ou,\{1,\ou\}^{2n+1}).
\end{eqnarray*}

\section*{Acknowledgment}

We thank the anonymous referee for suggestions which led to
improvements in the exposition.


\begin{thebibliography}{99}
\bibitem{AS}
Milton Abramowitz and Irene A.~Stegun, {\it Handbook of
Mathematical
Functions}, Dover, New York, 1972.

\bibitem{Bail}
David H.~Bailey, Jonathan M.~Borwein and Roland Girgensohn,
``Experimental Evaluation of Euler Sums,''
{\it Exp.\ Math.}, \textbf{3} (1994), 17--30.


\bibitem{BBG}
David Borwein, Jonathan M.~Borwein and Roland Girgensohn,
``Explicit Evaluation of Euler Sums,'' {\it Proc.\ Edinburgh\
Math.\ Soc.}, \textbf{38} (1995), 277--294.

\bibitem{BBB}
Jonathan M.~Borwein, David M.~Bradley and David J.~Broadhurst,
``Evaluations of $k$-fold Euler/Zagier Sums: A~Compendium of
Results for Arbitrary $k$,'' {\it Elec.\ J.\ Comb.}, \textbf{4}
(1997), no.~2, \#R5.

\bibitem{BBBLa}
Jonathan M.~Borwein, David M.~Bradley, David J.~Broadhurst and
Petr Lison\v ek, ``Special Values of Multiple Polylogarithms,''
{\it Trans.\ Amer.\ Math.\ Soc.} \textbf{355} (2001), no.~3,
907--941.

\bibitem{BBBLc}
Jonathan M.~Borwein, David M.~Bradley, David J.~Broadhurst and
Petr Lison\v ek, ``Combinatorial Aspects of Multiple Zeta
Values,'' {\it Elec.\ J.\ Comb.}, \textbf{5} (1998), no.~1, \#R38.

\bibitem{BG}
Jonathan M.~Borwein and Roland Girgensohn,
``Evaluation of Triple Euler Sums,''
{\it Elec.\ J.\ Comb.}, \textbf{3} (1996) \#R23.

\bibitem{BowBrad2}
Douglas Bowman and David M.~Bradley, ``The Algebra and
Combinatorics of Shuffles and Multiple Zeta Values'',  {\it J.\
Combinatorial Theory, Series A}, \textbf{97} (2002), 43--61.

\bibitem{BowBrad3}
Douglas Bowman and David M.~Bradley, ``Multiple Polylogarithms: A
Brief Survey,'' {\it Proceedings of a Conference on $q$-Series
with Applications to Combinatorics, Number Theory and Physics},
(Bruce C.~Berndt and Ken Ono eds.) American Mathematical Society,
Contemporary Mathematics \textbf{291} (2001), 71--92.

\bibitem{BowBradRyoo}
Douglas Bowm and, David M.~Bradley and Ji Hoon Ryoo, ``Some
Multi-Set Inclusions Associated with Shuffle Convolutions and
Multiple Zeta Values,'' {\it European J.\ Combinatorics},
\textbf{24} (2003), 121--127.

\bibitem{DJB1}
David J.~Broadhurst, ``On the Enumeration  of Irreducible $k$-fold
Euler Sums and Their Roles in Knot Theory and Field Theory,''
to appear in {\it J.\ Math.\ Phys.}

\bibitem{Bprivate}
David J.~Broadhurst, private e-mail, 1997.



\bibitem{BK1}
David J.\ Broadhurst and Dirk Kreimer, ``Association of Multiple
Zeta Values with Positive Knots via Feynman Diagrams up to 9
Loops,'' {\it Phys.\ Lett. B}, Vol.~393, 1997, pp.~403--412.









\bibitem{Drin}
V.\ G.\ Drinfeld, ``On Quasitriangular Quasi-Hopf Algebras and a
Group Closely Connected with Gal($\overline{\Q}/\Q$),'' {\it
Algebra i Analiz}, {\bf 2:4}(1990), 149--181. English transl.:
{\it Lenningrad Math.~J.}, Vol.~2, 1991, pp.~829--860.  [MR
92f:16047]


\bibitem{LE}
Leonhard Euler, ``Meditationes Circa Singulare Serierum Genus,''
{\it Novi Comm.\ Acad.\ Sci.\ Petropol.}, \textbf{20} (1775),
140--186, Reprinted in ``Opera Omnia'', ser.~I, vol.~15, B.
G.~Teubner, Berlin, 1927, pp.~217--267.






\bibitem{Gonch}
Alexander B.~Goncharov, ``Multiple Polylogarithms, Cyclotomy and
Modular Complexes,'' {\it Math.\ Res.\ Lett.}, \textbf{5} (1998),
no.~4, 497--516.

\bibitem{Gonch2}
Alexander B.~Goncharov, ``Polylogarithms in Arithmetic and
Geometry, '' Proc.\ ICM-94, Zurich, 1995, 374--387.



\bibitem{Hoff1}
Michael E.~Hoffman, ``Multiple Harmonic Series,'' {\it Pacific J.\
Math.}, \textbf{152} (1992), no.~2, 275--290.

\bibitem{Hoff2}
Michael E.~Hoffman, ``The Algebra of Multiple Harmonic Series,''
{\it J.\  Algebra}, \textbf{194} (1997), 477--495.

\bibitem{Hoff3}
Michael E.~Hoffman, ``Quasi-Shuffle Products,'' {\it J.\ Algebraic
Comb.}, \textbf{11} (2000), 49--68.

\bibitem{Hoff4}
Michael E.~Hoffman, and Yasuo Ohno, ``Relations of Multiple Zeta
Values and their Algebraic Expression,'' {\it J.\ Algebra},
\textbf{262} (2003), no.~2, 332--347.



\bibitem{Kass}
Christian Kassel, {\it Quantum Groups}, Springer-Verlag, New York,
1995.









\bibitem{MinPet1}
Hoang Ngoc Minh and Michel Petitot, ``Polylogarithms and the
Riemann $\z$ function,'' {\it Discrete Math.}, \textbf{217}
(2000), no.~1--3, 273--292.






\bibitem{YOhno}
Yasuo Ohno, ``A Generalization of the Duality and Sum Formulas on
the Multiple Zeta Values,'' {\it J.~Number Theory}, \textbf{74}
(1999), 39--43.

\bibitem{Rain}
Earl D.~Rainville, {\it Special Functions}, Chelsea Publishing,
New York, 1971.







\bibitem{Zag}
Don Zagier, ``Values of Zeta Functions and their Applications,''
{\it First European Congress of Mathematics}, Vol.~II,
Birkh\"auser, Boston, 1994, pp.~497--512.

\end{thebibliography}
\end{document}